# The Dantzig selector and sparsity oracle inequalities

VLADIMIR KOLTCHINSKII

*School of Mathematics, Georgia Institute of Technology, Atlanta, GA 30332-0160, USA.
E-mail: vlad@math.gatech.edu*

Let
$$Y_j = f_*(X_j) + \xi_j, \qquad j = 1, \ldots, n,$$
where $X, X_1, \ldots, X_n$ are i.i.d. random variables in a measurable space $(S, \mathcal{A})$ with distribution $\Pi$ and $\xi, \xi_1, \ldots, \xi_n$ are i.i.d. random variables with $\mathbb{E}\xi = 0$ independent of $(X_1, \ldots, X_n)$. Given a dictionary $h_1, \ldots, h_N : S \mapsto \mathbb{R}$, let $f_\lambda := \sum_{j=1}^N \lambda_j h_j$, $\lambda = (\lambda_1, \ldots, \lambda_N) \in \mathbb{R}^N$. Given $\varepsilon > 0$, define

$$\hat{\Lambda}_\varepsilon := \left\{ \lambda \in \mathbb{R}^N : \max_{1 \leq k \leq N} \left| n^{-1} \sum_{j=1}^n (f_\lambda(X_j) - Y_j) h_k(X_j) \right| \leq \varepsilon \right\}$$

and

$$\hat{\lambda} := \hat{\lambda}^\varepsilon \in \underset{\lambda \in \hat{\Lambda}_\varepsilon}{\operatorname{Argmin}} \|\lambda\|_{\ell_1}.$$

In the case where $f_* := f_{\lambda^*}, \lambda^* \in \mathbb{R}^N$, Candes and Tao [*Ann. Statist.* **35** (2007) 2313–2351] suggested using $\hat{\lambda}$ as an estimator of $\lambda^*$. They called this estimator "the Dantzig selector". We study the properties of $f_{\hat{\lambda}}$ as an estimator of $f_*$ for regression models with random design, extending some of the results of Candes and Tao (and providing alternative proofs of these results).

## 1. Introduction

Consider a regression model with random design,
$$Y_j = f_*(X_j) + \xi_j, \qquad j = 1, \ldots, n,$$
where $X, X_1, \ldots, X_n$ are i.i.d. random variables in a measurable space $(S, \mathcal{A})$ with distribution $\Pi$ and $\xi, \xi_1, \ldots, \xi_n$ are i.i.d. random variables with $\mathbb{E}\xi = 0$, independent of $(X_1, \ldots, X_n)$ (in what follows, it will be assumed that the noise $\xi_j$ satisfies some further assumptions, such as, for instance, $\xi_j$ is $N(0, \sigma^2)$).

Let $h_1, \ldots, h_N$ be a dictionary consisting of $N \geq 2$ functions from $S$ into $\mathbb{R}$. Define

$$f_\lambda := \sum_{j=1}^{N} \lambda_j h_j, \qquad \lambda = (\lambda_1, \ldots, \lambda_N) \in \mathbb{R}^N.$$

Given $\varepsilon > 0$, define the set

$$\hat{\Lambda}_\varepsilon := \left\{ \lambda \in \mathbb{R}^N : \max_{1 \leq k \leq N} \left| n^{-1} \sum_{j=1}^{n} (f_\lambda(X_j) - Y_j) h_k(X_j) \right| \leq \varepsilon \right\}$$

and consider

$$\hat{\lambda} := \hat{\lambda}^\varepsilon \in \operatorname*{Argmin}_{\lambda \in \hat{\Lambda}_\varepsilon} \|\lambda\|_{\ell_1}.$$

Although the set of constraints $\hat{\Lambda}_\varepsilon$ could be empty, we will see that for sufficiently large values of $\varepsilon$, it is non-empty with a high probability (if $\hat{\Lambda}_\varepsilon = \varnothing$, one can define $\hat{\lambda}^\varepsilon$ in an arbitrary way, for instance, $\hat{\lambda}^\varepsilon = 0$).

In the case where $f_* = f_{\lambda^*}$ for some $\lambda^* \in \mathbb{R}^N$, Candes and Tao (2007) suggested using $\hat{\lambda}^\varepsilon$ as an estimator of the vector of coefficients $\lambda^*$. It is easy to see that the computation of $\hat{\lambda}^\varepsilon$ reduces to a linear programming problem:

$$\sum_{j=1}^{N} u_j \to \min,$$

subject to the constraints

$$u_k \geq 0, \qquad -u_k \leq \lambda_k \leq u_k, \qquad -\varepsilon \leq n^{-1} \sum_{j=1}^{n} (f_\lambda(X_j) - Y_j) h_k(X_j) \leq \varepsilon, \qquad k = 1, \ldots, N.$$

Candes and Tao called this estimator "the Dantzig selector". It is closely related to the $\ell_1$-penalization method (similar to what is called "LASSO" in statistical literature), which is based on fitting the regression model by solving the following penalized empirical risk minimization problem:

$$n^{-1} \sum_{j=1}^{n} (f_\lambda(X_j) - Y_j)^2 + 2\varepsilon \|\lambda\|_{\ell_1} =: L_n(\lambda) + 2\varepsilon \|\lambda\|_{\ell_1} \to \min. \qquad (1.1)$$

Note that

$$\hat{\Lambda}_\varepsilon = \{\lambda : \|\nabla L_n(\lambda)\|_{\ell_\infty} \leq 2\varepsilon\}$$

and that $\lambda \in \hat{\Lambda}_\varepsilon$ is a necessary condition for $\lambda$ to be a solution of (1.1).

We will establish several "sparsity oracle inequalities" for the Dantzig selector that are akin to recent inequalities proved in Bunea, Tsybakov and Wegkamp (2007), van de Geer (2008) and Koltchinskii (2009) in the case of $\ell_1$- or $\ell_p$-penalized empirical risk minimization. Candes and Tao (2007) concentrated on the case of fixed design regression models, that is, when the design points $X_1, \ldots, X_n$ are non-random. They proved their version of oracle inequalities under the basic assumption that the design matrix $A = (h_j(X_i))_{i=1,n;j=1,N}$ satisfies the so called uniform uncertainty principle (UUP). To explain the meaning of this assumption, define

$$J_\lambda := \mathrm{supp}(\lambda) := \{j : \lambda_j \neq 0\}, \qquad \lambda \in \mathbb{R}^N,$$

and set $d(\lambda) := \mathrm{card}(J_\lambda)$. Define $\delta_d(\Pi)$ to be the smallest $\delta > 0$ such that for all $\lambda \in \mathbb{R}^N$ with $d(\lambda) \leq d$,

$$(1-\delta)\|\lambda\|_{\ell_2} \leq \left\|\sum_{j=1}^N \lambda_j h_j\right\|_{L_2(\Pi)} \leq (1+\delta)\|\lambda\|_{\ell_2}.$$

If $\delta_d(\Pi) < 1$, then $d$-dimensional subspaces spanned on subsets of the dictionary and equipped with either the $L_2(\Pi)$-norm, or the $\ell_2$-norm on vectors of coefficients are "almost" isometric. Given the dictionary $\{h_1, \ldots, h_N\}$, it is natural to call the quantity $\delta_d(\Pi)$ *the restricted isometry constant* of dimension $d$ with respect to measure $\Pi$. If $\Pi_n$ denotes the empirical measure based on the design points $X_1, \ldots, X_n$, then the UUP essentially means that the restricted isometry constants $\delta_d(\Pi_n)$ (which are characteristics of the design matrix $A$) are sufficiently small for the values of $d$ comparable with the degree of sparsity of representation of $f_*$ in the dictionary (the number of non-zero coefficients of $\lambda^*$). Candes and Tao (2007) stated that the UUP holds with a high probability for some random design matrices such as the *Gaussian ensemble* (the matrix with i.i.d. standard normal entries). It is also true for the *Bernoulli* or *Rademacher ensemble* (the matrix with i.i.d. entries taking values $+1$ and $-1$ with probability $1/2$), which relies on some facts concerning random matrices that were established in other papers. In these examples, the dictionaries are orthonormal systems in the space $L_2(\Pi)$, which means that $\delta_d(\Pi) = 0$.

We here provide more direct proofs of oracle inequalities in the random design case that do not rely on the bounds for random matrices and that apply to broader classes of design distributions, in particular, to such distributions that the dictionary is not necessarily orthonormal in $L_2(\Pi)$, but rather satisfies a restricted isometry condition with respect to $\Pi$. The next statement is a typical example of what follows from the results of Sections 2 and 3 (specifically, from Corollary 6).

**Proposition 1.** *Suppose that the random vector* $(h_1(X), \ldots, h_N(X))$ *has normal distribution with zero mean and that the noise $\xi$ is $N(0; \sigma^2)$. In addition, suppose that $f_* = f_{\lambda^*}, \lambda^* \in \mathbb{R}^N$. There then exist constants $\bar{\delta} \in (0,1)$ and $C, D > 0$ with the following property. For an arbitrary $A \geq 1$, denote by $\bar{d}$ the largest $d \leq N/\mathrm{e} - 1$ such that*

$$\delta_{3d}(\Pi) \leq \bar{\delta}$$

*and*

$$C\sqrt{\frac{Ad\log(N/d)}{n}} \le 1/4.$$

*Then, for all*

$$\varepsilon \ge C\sigma \max_{1\le k\le N} \|h_k\|_{L_2(\Pi)} \sqrt{\frac{A\log N}{n}},$$

*the condition $d(\lambda^*) \le \bar{d}$ implies that with probability at least $1 - N^{-A}$,*

$$\|\hat{\lambda} - \lambda^*\|_{\ell_2} \le D\sqrt{d(\lambda^*)}\varepsilon.$$

Our approach is based on some facts concerning empirical and Rademacher processes and it is close to the approach taken by Rudelson and Vershynin (2005) or Mendelson, Pajor and Tomczak-Jaegermann (2007). At the same time, it relies only on rather elementary tools (symmetrization and contraction inequalities for Rademacher processes and Bernstein-type exponential bounds) and does not use more advanced techniques, such as concentration of measure and generic chaining, which are used in the papers cited above. It is worth mentioning that Koltchinskii (2005, 2009) showed that if, in (1.1), one uses $\|\lambda\|_{\ell_p}^p$ with $p = 1 + \frac{c}{\log N}$ instead of $\|\lambda\|_{\ell_1}$, then one can establish a version of sparsity oracle inequalities without making strong assumptions on the dictionary such as a restricted isometry condition.

In the next section, we introduce some geometric characteristics of the dictionary that are of importance in analysis of sparse recovery problems, and we prove general oracle inequalities for the Dantzig selector in terms of these characteristics. Several corollaries, more special results and some examples are given in Section 3. Finally, the Appendix contains some exponential bounds for Rademacher processes needed in the proofs of the main results.

## 2. Main results

In what follows, we frequently use the Orlicz norm $\|\cdot\|_\psi$ for random variables, most often with $\psi = \psi_1$, $\psi_1(x) := e^{|x|} - 1$ or $\psi = \psi_2$, $\psi_2(x) := e^{x^2} - 1$. For any convex non-decreasing function $\psi : \mathbb{R}_+ \mapsto \mathbb{R}_+$ with $\psi(0) = 0$, it is defined as

$$\|\eta\|_\psi := \inf\left\{C > 0 : \mathbb{E}\psi\left(\frac{|\eta|}{C}\right) \le 1\right\}$$

(see Ledoux and Talagrand (1991), van der Vaart and Wellner (1996), de la Pena and Giné (1998)).

For $J \subset \{1, \ldots, N\}$, let $d(J) := \text{card}(J)$. Define

$$C_J := \left\{u \in \mathbb{R}^N : \sum_{j \notin J} |u_j| \le \sum_{j \in J} |u_j|\right\}.$$

The set $C_J$ is a cone in $\mathbb{R}^N$ (that is, $u \in C_J$ implies that $\alpha u \in C_J$ for all $\alpha \geq 0$). It consists of vectors $u \in \mathbb{R}^N$ such that the coordinates of $u$ in the set $J$ are *dominant*. Such *cones of dominant coordinates* play an important role in the analysis of the Dantzig selector, LASSO and other sparse recovery methods. The reason is that, for a "sparse" feasible vector $\lambda \in \hat{\Lambda}_\varepsilon$, the definition of the Dantzig selector $\hat{\lambda}^\varepsilon$ means that $\|\hat{\lambda}^\varepsilon\|_{\ell_1} \leq \|\lambda\|_{\ell_1}$, which implies that

$$\sum_{j \notin J_\lambda} |\hat{\lambda}^\varepsilon_j - \lambda_j| = \sum_{j \notin J_\lambda} |\hat{\lambda}^\varepsilon_j| \leq \sum_{j \in J_\lambda} (|\lambda_j| - |\hat{\lambda}^\varepsilon_j|) \leq \sum_{j \in J_\lambda} |\hat{\lambda}^\varepsilon_j - \lambda_j| \qquad (2.1)$$

and, hence, $\hat{\lambda}^\varepsilon - \lambda \in C_{J_\lambda}$. The proofs of various bounds on the norms of the vector $\hat{\lambda}^\varepsilon - \lambda$ and on the norms of the corresponding function

$$f_{\hat{\lambda}^\varepsilon} - f_\lambda \in \text{l.s.}(\{h_1, \ldots, h_N\})$$

are usually based on the comparison of these norms on the cone $C_{J_\lambda}$. We will introduce several geometric characteristics of the dictionary that are needed for such a comparison.

Define

$$\beta(J) := \beta(J; \Pi) := \inf\left\{\beta > 0 : \forall \lambda \in C_J, \sum_{j \in J} |\lambda_j| \leq \beta \left\|\sum_{j=1}^N \lambda_j h_j\right\|_{L_1(\Pi)}\right\}$$

(set $\beta(J) = 0$ if $J = \varnothing$). Note that if $J \neq \varnothing$ and the functions $h_1, \ldots, h_N$ in the dictionary are linearly independent in $L_1(\Pi)$, then $\beta(J) < +\infty$.

Another quantity of interest is

$$\beta_2(J) := \beta_2(J; \Pi) := \inf\left\{\beta > 0 : \forall \lambda \in C_J, \sum_{j \in J} |\lambda_j|^2 \leq \beta^2 \left\|\sum_{j=1}^N \lambda_j h_j\right\|_{L_2(\Pi)}^2\right\}.$$

Note that $\beta_2(J) = 1, J \neq \varnothing$ if the dictionary $\{h_1, \ldots, h_N\}$ is orthonormal. In Section 3, the connection of these quantities to restricted isometry constants $\delta_d(\Pi)$ is discussed. In particular, it can be shown that if $\delta_{3d}(\Pi)$ is small enough, then, for all sets $J$ of cardinality $d$, $\beta_2(J)$ remains properly bounded.

The following condition on the dictionary and on the distribution $\Pi$ is often of interest: for all $\lambda \in C_J$

$$\left\|\sum_{j=1}^N \lambda_j h_j\right\|_{L_1(\Pi)} \leq \left\|\sum_{j=1}^N \lambda_j h_j\right\|_{L_2(\Pi)} \leq B(J) \left\|\sum_{j=1}^N \lambda_j h_j\right\|_{L_1(\Pi)} \qquad (2.2)$$

with some constant $B(J) > 0$.

Note that the first inequality in (2.2) is trivial for all $\lambda \in \mathbb{R}^N$. The second (non-trivial) bound holds for all $\lambda \in \mathbb{R}^N$, with a constant $B > 0$ that does not depend on $N$, and on the set $J$ in several interesting, but rather special, examples. In particular, this condition holds when $(h_1(X), \ldots, h_N(X))$ has mean zero normal distribution in $\mathbb{R}^N$ (for instance, if $h_1(X), \ldots, h_N(X)$ are i.i.d. standard normal, which is the case for the Gaussian dictionary) or when $h_1(X), \ldots, h_N(X)$ are i.i.d. Rademacher random variables (that is, $h_j(X)$ is $+1$ or $-1$ with probability $1/2$ each; this is the case for Bernoulli and Rademacher dictionaries). In the last case, (2.2) holds by the Khinchine inequality. For Gaussian and Bernoulli dictionaries, all $L_p$-norms, $p \geq 1$, and even $\psi_1$- and $\psi_2$-norms of $\sum_{j=1}^{N} \lambda_j h_j$ are equivalent up to numerical constants (see Bobkov and Houdré (1997) for a discussion of more general Khinchine-type inequalities and their connections with isoperimetric constants). In general, the constant $B$ might depend on $N$ and, since this constant is involved in the bounds on the performance of the Dantzig selector, it is of some importance that condition (2.2) is supposed to hold only for all $\lambda \in C_J$ (rather than for all $\lambda \in \mathbb{R}^N$), and this condition is usually needed for a small set $J$.

Under the condition (2.2), the following bound is straightforward:

$$\beta(J) \leq B(J)\beta_2(J)\sqrt{d(J)}. \qquad (2.3)$$

If $\beta_2(J)$ is bounded by a small constant (as in the case of orthonormal dictionaries), then $\beta(J)$ is "small" for sets $J$ of small cardinality $d(J)$.

Recall the notation $J_\lambda := \mathrm{supp}(\lambda)$ and also recall that $d(\lambda) := d(J_\lambda)$.

We will fix the values of $\varepsilon > 0$, $A > 0$ and $C > 0$, assume that

$$\frac{A \log N}{n} \leq 1$$

and define the following set:

$$\Lambda_\varepsilon(A) := \left\{ \lambda \in \mathbb{R}^N : |\langle f_\lambda - f_*, h_k \rangle_{L_2(\Pi)}| \right.$$

$$\left. + C(\|(f_\lambda - f_*)(X) h_k(X)\|_{\psi_1} + \|\xi h_k(X)\|_{\psi_1}) \sqrt{\frac{A \log N}{n}} \leq \varepsilon, \ k = 1, \ldots, N \right\}$$

(recall that $\xi$ involved in the above definition is the noise of the regression model). Under the condition

$$\varepsilon \geq C \max_{1 \leq k \leq N} \|\xi h_k(X)\|_{\psi_1} \sqrt{\frac{A \log N}{n}},$$

which is necessary for the set $\Lambda_\varepsilon(A)$ to be non-empty, this set consists of vectors $\lambda$ such that $f_\lambda$ is, in a certain sense, a good approximation of $f_*$. The condition

$$\max_{1 \leq k \leq N} |\langle f_\lambda - f_*, h_k \rangle_{L_2(\Pi)}| \leq \varepsilon \qquad (2.4)$$

that follows from $\lambda \in \Lambda_\varepsilon(A)$ essentially means that $f_\lambda - f_*$ is almost ("up to $\varepsilon$") orthogonal to the linear span of the dictionary, so $f_\lambda$ should be close to the projection of $f_*$ on the linear span. In fact, (2.4) is a necessary condition of the minimum in the convex minimization problem

$$\|f_\lambda - f_*\|^2_{L_2(\Pi)} + 2\varepsilon\|\lambda\|_{\ell_1} \to \min,$$

which can be viewed as a distribution-dependent version of the empirical risk minimization (1.1) (recall that $\lambda \in \hat{\Lambda}_\varepsilon$ is a necessary condition for (1.1)).

If $f_{\lambda^0}, \lambda^0 \in \mathbb{R}^N$ is the orthogonal projection in $L_2(\Pi)$ of the function $f_*$ onto the linear span of the dictionary, then it is obvious that the condition

$$\varepsilon \geq C \max_{1\leq k\leq N}(\|(f_{\lambda^0} - f_*)(X)h_k(X)\|_{\psi_1} + \|\xi h_k(X)\|_{\psi_1})\sqrt{\frac{A\log N}{n}}$$

is sufficient for $\Lambda_\varepsilon(A) \neq \varnothing$ (since, under this condition, $\lambda^0 \in \Lambda_\varepsilon(A)$).

The next proposition shows that if $\Lambda_\varepsilon(A) \neq \varnothing$, then, with a high probability, $\hat{\Lambda}_\varepsilon \neq \varnothing$.

**Proposition 2.** *Suppose that $\Lambda_\varepsilon(A) \neq \varnothing$ and $\lambda \in \Lambda_\varepsilon(A)$. Then, with probability at least $1 - 2N^{-A}$, $\lambda \in \hat{\Lambda}_\varepsilon$.*

**Proof.** Indeed, for any such $\lambda$, we have

$$\left|n^{-1}\sum_{j=1}^n (f_\lambda(X_j) - Y_j)h_k(X_j)\right|$$

$$\leq |\langle f_\lambda - f_*, h_k\rangle_{L_2(\Pi)}|$$

$$+ \left|n^{-1}\sum_{j=1}^n [(f_\lambda(X_j) - f_*(X_j))h_k(X_j) - \mathbb{E}(f_\lambda(X) - f_*(X))h_k(X)]\right|$$

$$+ \left|n^{-1}\sum_{j=1}^n \xi_j h_k(X_j)\right|.$$

Applying Lemma 3 from the Appendix to the second and third terms yields, with probability at least $1 - 2N^{-A}$,

$$\max_{1\leq k\leq N}\left|n^{-1}\sum_{j=1}^n (f_\lambda(X_j) - Y_j)h_k(X_j)\right|$$

$$\leq \max_{1\leq k\leq N}\left[|\langle f_\lambda - f_*, h_k\rangle_{L_2(\Pi)}| + C(\|(f_\lambda - f_*)(X)h_k(X)\|_{\psi_1} + \|\xi h_k(X)\|_{\psi_1})\sqrt{\frac{A\log N}{n}}\right]$$

$$\leq \varepsilon,$$

by the definition of the set $\Lambda_\varepsilon(A)$. □

Also, define

$$\Lambda_S(A) := \left\{ \lambda \in \mathbb{R}^N : C\beta(J_\lambda) \max_{1 \leq k \leq N} \|h_k(X)\|_{\psi_1} \sqrt{\frac{A \log N}{n}} \leq 1/4 \right\}.$$

We will interpret $\Lambda_S(A)$ as a set of "sparse" vectors since, in view of the bound (2.3), $\beta(J_\lambda)$ has some connection to the sparsity of $\lambda$. Of course, the fact that $\beta(J_\lambda)$ is not too large is also related to the properties of the dictionary. For dictionaries that are close to being orthonormal, $\Lambda_S(A)$ would include sparse enough vectors in the usual sense.

Essentially, the bounds of Theorems 1 and 2 below show that if there exists a vector $\lambda$ in $\hat{\Lambda}_\varepsilon$ (the set of constraints of the Dantzig selector) that is sufficiently "sparse", then the Dantzig selector will be in a small ball around $\lambda$ in such norms as $\|\cdot\|_{\ell_1}$ and $\|\cdot\|_{\ell_2}$, or $f_{\hat{\lambda}}$ will be in a small ball around $f_\lambda$ in such norms as $\|\cdot\|_{L_1(\Pi)}$ and $\|\cdot\|_{L_2(\Pi)}$. The radius of this ball crucially depends on the degree of sparsity of $\lambda$ and also on the "well-posedness" of the dictionary characterized by such quantities as $\beta_2$ (see also Section 3 for a discussion of the connection of these quantities to restricted isometry properties of the dictionary). The bounds also imply that the Dantzig selector is adaptive to an unknown degree of the sparsity of the problem (at least in the case when the dictionary is not very far from being orthonormal in $L_2(\Pi)$).

Let

$$\tilde{\Lambda}_\varepsilon(A) := \Lambda_\varepsilon(A) \cap \Lambda_S(A).$$

This set will be interpreted in the next theorem as a set of oracle vectors and it will be assumed that $\tilde{\Lambda}_\varepsilon(A) \neq \varnothing$. In particular, it means that $\varepsilon$ must satisfy

$$\varepsilon \geq C \max_{1 \leq k \leq N} \|\xi h_k(X)\|_{\psi_1} \sqrt{\frac{A \log N}{n}}$$

(which, of course, requires that $\|\xi h_k(X)\|_{\psi_1} < +\infty$). The fact that $\lambda \in \tilde{\Lambda}_\varepsilon(A)$ implies that $\lambda$ is sparse in the sense that $\lambda \in \Lambda_S(A)$ and, at the same time, that $f_\lambda$ provides a reasonably good approximation of $f_*$ in the sense that both (2.4) holds and

$$C\|(f_\lambda - f_*)(X)h_k(X)\|_{\psi_1} \sqrt{\frac{A \log N}{n}} \leq \varepsilon.$$

If $\tilde{\Lambda}_\varepsilon(A) = \varnothing$, then there are no sparse vectors $\lambda$ for which $f_\lambda$ approximates $f_*$ well, so, from this point of view, the problem is not sparse.

First, we prove the following result.

**Theorem 1.** *There exists a constant $C$ in the definitions of $\Lambda_\varepsilon(A), \Lambda_S(A)$ such that for all $A \geq 1$ with probability at least $1 - N^{-A}$, the following bounds hold for all $\lambda \in \hat{\Lambda}_\varepsilon \cap \Lambda_S(A)$ and for the Dantzig selector $\hat{\lambda}$:*

$$\|f_{\hat{\lambda}} - f_\lambda\|_{L_1(\Pi)} \leq 16\beta(J_\lambda)\varepsilon$$

and
$$\|\hat{\lambda} - \lambda\|_{\ell_1} \leq 32\beta^2(J_\lambda)\varepsilon.$$

Under the assumption that $\tilde{\Lambda}_\varepsilon(A) \neq 0$, with the same probability,
$$\|f_{\hat{\lambda}} - f_*\|_{L_1(\Pi)} \leq \inf_{\lambda \in \tilde{\Lambda}_\varepsilon(A)} [\|f_\lambda - f_*\|_{L_1(\Pi)} + 16\beta(J_\lambda)\varepsilon],$$

and if, in addition, $f_* = f_{\lambda^*}, \lambda^* \in \mathbb{R}^N$, then we also have
$$\|\hat{\lambda} - \lambda^*\|_{\ell_1} \leq \inf_{\lambda \in \tilde{\Lambda}_\varepsilon(A)} [\|\lambda - \lambda^*\|_{\ell_1} + 32\beta^2(J_\lambda)\varepsilon].$$

**Proof.** Suppose that $\lambda \in \hat{\Lambda}_\varepsilon \cap \Lambda_S(A)$. The proof of the first two bounds will be based on upper bounding $\|\hat{\lambda} - \lambda\|_{\ell_1}$ in terms of $\|f_{\hat{\lambda}} - f_\lambda\|_{L_1(\Pi)}$ and vice versa. Combining these bounds yields inequalities on both of the norms that can be solved, leading to the first two bounds of the theorem.

Since $\lambda \in \hat{\Lambda}^\varepsilon$, (2.1) implies that $\hat{\lambda} - \lambda \in C_{J_\lambda}$ and

$$\|\hat{\lambda} - \lambda\|_{\ell_1} \leq \sum_{j \notin J_\lambda} |\hat{\lambda}_j| + \sum_{j \in J_\lambda} |\lambda_j - \hat{\lambda}_j| \leq 2 \sum_{j \in J_\lambda} |\lambda_j - \hat{\lambda}_j| \leq 2\beta(J_\lambda)\|f_{\hat{\lambda}} - f_\lambda\|_{L_1(\Pi)}. \quad (2.5)$$

We will now upper bound $\|f_{\hat{\lambda}} - f_\lambda\|_{L_1(\Pi)}$ in terms of $\|\hat{\lambda} - \lambda\|_{\ell_1}$, which will imply the result. We start with the following, obvious, bound (we use the notation $\nu(f) := \int f \, d\nu$):

$$\|f_{\hat{\lambda}} - f_\lambda\|_{L_1(\Pi)} = \|f_{\hat{\lambda}} - f_\lambda\|_{L_1(\Pi_n)} + (\Pi - \Pi_n)(|f_{\hat{\lambda}} - f_\lambda|)$$
$$\leq \|f_{\hat{\lambda}} - f_\lambda\|_{L_1(\Pi_n)} + \sup_{\|u\|_{\ell_1} \leq 1} |(\Pi_n - \Pi)(|f_u|)| \|\hat{\lambda} - \lambda\|_{\ell_1}. \quad (2.6)$$

We will separately bound the first and second terms of this bound. First, note that

$$\|f_{\hat{\lambda}} - f_\lambda\|_{L_1(\Pi_n)}^2 \leq \|f_{\hat{\lambda}} - f_\lambda\|_{L_2(\Pi_n)}^2 = \langle f_{\hat{\lambda}} - f_\lambda, f_{\hat{\lambda}} - f_\lambda \rangle_{L_2(\Pi_n)}$$
$$= \sum_{k=1}^N (\hat{\lambda}_k - \lambda_k) \langle f_{\hat{\lambda}} - f_\lambda, h_k \rangle_{L_2(\Pi_n)} \leq \|\hat{\lambda} - \lambda\|_{\ell_1} \max_{1 \leq k \leq N} |\langle f_{\hat{\lambda}} - f_\lambda, h_k \rangle_{L_2(\Pi_n)}|.$$

Since both $\hat{\lambda} \in \hat{\Lambda}_\varepsilon$ and $\lambda \in \hat{\Lambda}_\varepsilon$, we have

$$\max_{1 \leq k \leq N} |\langle f_{\hat{\lambda}} - f_\lambda, h_k \rangle_{L_2(\Pi_n)}|$$
$$\leq \max_{1 \leq k \leq N} \left| n^{-1} \sum_{j=1}^n (f_\lambda(X_j) - Y_j) h_k(X_j) \right| + \max_{1 \leq k \leq N} \left| n^{-1} \sum_{j=1}^n (f_{\hat{\lambda}}(X_j) - Y_j) h_k(X_j) \right| \leq 2\varepsilon,$$

which implies that
$$\|f_{\hat{\lambda}} - f_\lambda\|_{L_1(\Pi_n)} \leq \sqrt{2\varepsilon \|\hat{\lambda} - \lambda\|_{\ell_1}}.$$

By Lemma 4 from the Appendix, with probability at least $1 - N^{-A}$ (under the assumption $A \log N \leq n$),

$$\sup_{\|u\|_{\ell_1} \leq 1} |(\Pi_n - \Pi)(|f_u|)| \leq C \max_{1 \leq k \leq N} \|h_k\|_{\psi_1} \sqrt{\frac{A \log N}{n}}.$$

This yields the following bound (with probability at least $1 - N^{-A}$):

$$\|f_{\hat{\lambda}} - f_{\lambda}\|_{L_1(\Pi)} \leq \sqrt{2\varepsilon \|\hat{\lambda} - \lambda\|_{\ell_1}} + C \max_{1 \leq k \leq N} \|h_k\|_{\psi_1} \sqrt{\frac{A \log N}{n}} \|\hat{\lambda} - \lambda\|_{\ell_1}. \quad (2.7)$$

Together with (2.5), this implies that

$$\|f_{\hat{\lambda}} - f_{\lambda}\|_{L_1(\Pi)} \leq \sqrt{4\varepsilon \beta(J_{\lambda}) \|f_{\hat{\lambda}} - f_{\lambda}\|_{L_1(\Pi)}}$$
$$+ 2C \max_{1 \leq k \leq N} \|h_k\|_{\psi_1} \sqrt{\frac{A \log N}{n}} \beta(J_{\lambda}) \|f_{\hat{\lambda}} - f_{\lambda}\|_{L_1(\Pi)}.$$

Recalling the definition of the set $\Lambda_S(A)$, we can guarantee that

$$2C \max_{1 \leq k \leq N} \|h_k\|_{\psi_1} \sqrt{\frac{A \log N}{n}} \beta(J_{\lambda}) \leq 1/2,$$

which implies that

$$\|f_{\hat{\lambda}} - f_{\lambda}\|_{L_1(\Pi)} \leq 2\sqrt{4\varepsilon \beta(J_{\lambda}) \|f_{\hat{\lambda}} - f_{\lambda}\|_{L_1(\Pi)}},$$

and the first bound now follows. The second bound is also true, in view of (2.5).

To prove each of the remaining bounds, define $\bar{\lambda}$ to be the vector for which the infimum in the right-hand side of the bound is attained. By Proposition 2, with probability at least $1 - 2N^{-A}$, we have $\bar{\lambda} \in \hat{\Lambda}_{\varepsilon} \cap \Lambda_S(A)$. Therefore, we can use the first two bounds of the theorem and the triangle inequality to complete the proof of the remaining bounds that now hold with probability at least $1 - 3N^{-A}$.

It only remains to show that by adjusting the value of the constant $C$ in the definitions of the sets $\Lambda_{\varepsilon}(A)$ and $\Lambda_S(A)$, it is possible to ensure that the bounds hold with probability at least $1 - N^{-A}$, as was claimed. To this end, check that for $c := \log_2 3 + 1$ and all $A \geq c, N \geq 2$,

$$3N^{-A} \leq N^{-A/c}.$$

Now, take $A' = A/c \geq 1$ and replace the constant $C$ with $C\sqrt{c}$ to show that the bounds hold with probability at least $1 - N^{-A'}$. $\square$

Under the condition (2.2), the bound (2.3) holds and one can derive from Theorem 1 the bounds expressed in terms of quantity $\beta_2$, namely, replacing $\beta(J_{\lambda})$ in the inequalities of Theorem 1 by the upper bound $B(J_{\lambda})\beta_2(J_{\lambda})\sqrt{d(\lambda)}$. However, below, we will give another version of such a statement with bounds on the norms $\|\cdot\|_{L_2(\Pi)}$ and $\|\cdot\|_{\ell_2}$, and with a slight improvement of the logarithmic factor in the definition of the set of sparse vectors $\Lambda_S(A)$.

Define (with a minor abuse of notation)

$$\beta_2(d) := \beta_2(d;\Pi) := \max\{\beta_2(J) : J \subset \{1,\ldots,N\}, d(J) \leq 2d\}$$

and

$$B(d) := \max\{B(J) : J \subset \{1,\ldots,N\}, d(J) \leq d\}.$$

Let $\bar{d}$ denote the largest $d \leq \frac{N}{e} - 1$ such that

$$\frac{A d \log(N/d)}{n} \leq 1$$

and

$$CB(d)\beta_2(d) \sup_{\|u\|_{\ell_2} \leq 1, d(u) \leq d} \|f_u\|_{\psi_1} \sqrt{\frac{A d \log(N/d)}{n}} \leq 1/4.$$

We redefine the set of "sparse" vectors as follows:

$$\Lambda_{S,2}(A) := \{\lambda \in \mathbb{R}^N : d(\lambda) \leq \bar{d}\}.$$

Let

$$\tilde{\Lambda}_\varepsilon^2(A) := \Lambda_\varepsilon(A) \cap \Lambda_{S,2}(A),$$

which will now play the role of an oracle set (that is, the set of sparse enough vectors that approximate the target function $f_*$ reasonably well).

We will use the notation

$$\binom{n}{\leq k} = \sum_{j=0}^{k} \binom{n}{j}.$$

**Theorem 2.** *Suppose that condition (2.2) holds. There exists a constant $C$ in the definitions of $\Lambda_\varepsilon(A)$ and $\Lambda_{S,2}(A)$ such that for all $A \geq 1$, with probability at least*

$$1 - 5^{-\bar{d}A}\binom{N}{\leq \bar{d}}^{-A},$$

*the following bounds hold:* $\forall \lambda \in \hat{\Lambda}_\varepsilon \cap \Lambda_{S,2}(A)$

$$\|f_{\hat{\lambda}} - f_\lambda\|_{L_2(\Pi)} \leq 16 B^2(d(\lambda))\beta_2(d(\lambda))\sqrt{d(\lambda)}\varepsilon$$

*and*

$$\|\hat{\lambda} - \lambda\|_{\ell_2} \leq 32 B^2(d(\lambda))\beta_2^2(d(\lambda))\sqrt{d(\lambda)}\varepsilon.$$

*Suppose that $\tilde{\Lambda}_\varepsilon^2(A) \neq \varnothing$. Then, with probability at least $1 - N^{-A}$,*

$$\|f_{\hat{\lambda}} - f_*\|_{L_2(\Pi)} \leq \inf_{\lambda \in \tilde{\Lambda}_\varepsilon^2(A)} [\|f_\lambda - f_*\|_{L_2(\Pi)} + 16 B^2(d(\lambda)) \beta_2(d(\lambda)) \sqrt{d(\lambda)} \varepsilon].$$

*Moreover, if, in addition, $f_* = f_{\lambda^*}, \lambda^* \in \mathbb{R}^N$, then, also with the same probability,*

$$\|\hat{\lambda} - \lambda^*\|_{\ell_2} \leq \inf_{\lambda \in \tilde{\Lambda}_\varepsilon^2(A)} [\|\lambda - \lambda^*\|_{\ell_2} + 32 B^2(d(\lambda)) \beta_2^2(d(\lambda)) \sqrt{d(\lambda)} \varepsilon].$$

We will need the following well-known fact (see Candes and Tao (2005), proof of their Theorem 1).

**Lemma 1.** *Let $u \in C_J$. Define $J_0 := J$ and let $J_1$ be the set of $d$ coordinates in $\{1, \ldots, N\} \setminus J_0$ for which the $|u_j|$'s are the largest, $J_2$ be the set of $d$ coordinates in $\{1, \ldots, N\} \setminus (J_0 \cup J_1)$ for which the $|u_j|$'s are the largest, etcetera. Define $u^{(k)} := (u_j : j \in J_k)$. Then, $u = \sum_{k \geq 0} u^{(k)}$ and*

$$\sum_{k \geq 2} \|u^{(k)}\|_{\ell_2} \leq \left( \sum_{j \in J} |u_j|^2 \right)^{1/2},$$

*which also implies that*

$$\|u\|_{\ell_2} \leq 2 \left( \sum_{j \in J_0 \cup J_1} |u_j|^2 \right)^{1/2}.$$

**Proof.** For all $j \in J_{k+1}$,

$$|u_j| \leq \frac{1}{d} \sum_{i \in J_k} |u_i|,$$

implying that

$$\left( \sum_{j \in J_{k+1}} |u_j|^2 \right)^{1/2} \leq \frac{1}{\sqrt{d}} \sum_{j \in J_k} |u_j|.$$

Adding these inequalities for $k = 1, 2, \ldots$ yields

$$\sum_{k \geq 2} \|u^{(k)}\|_{\ell_2} \leq \frac{1}{\sqrt{d}} \sum_{j \notin J} |u_j| \leq \frac{1}{\sqrt{d}} \sum_{j \in J} |u_j| \leq \left( \sum_{j \in J} |u_j|^2 \right)^{1/2} \leq \left( \sum_{j \in J \cup J_1} |u_j|^2 \right)^{1/2}.$$

Thus, for $u \in C_J$,

$$\|u\|_{\ell_2} \leq 2 \left( \sum_{j \in J_0 \cup J_1} |u_j|^2 \right)^{1/2}.$$

□

**Proof of Theorem 2.** This is a straightforward modification of the proof of Theorem 1. Suppose that $\lambda \in \hat{\Lambda}_\varepsilon \cap \Lambda_{S,2}(A)$. Instead of (2.6), we now use

$$\begin{aligned}\|f_{\hat{\lambda}} - f_\lambda\|_{L_1(\Pi)} &= \|f_{\hat{\lambda}} - f_\lambda\|_{L_1(\Pi_n)} + (\Pi - \Pi_n)(|f_{\hat{\lambda}} - f_\lambda|) \\ &\leq \|f_{\hat{\lambda}} - f_\lambda\|_{L_1(\Pi_n)} + \sup_{\|u\|_{\ell_2} \leq 1, u \in C_{J_\lambda}} |(\Pi_n - \Pi)(|f_u|)| \|\hat{\lambda} - \lambda\|_{\ell_2}.\end{aligned} \quad (2.8)$$

To bound $\|\hat{\lambda} - \lambda\|_{\ell_2}$, we observe (as in the proof of Theorem 1) that $\hat{\lambda} - \lambda \in C_{J_\lambda}$ and apply Lemma 1 to $u = \hat{\lambda} - \lambda$, $J = J_\lambda$:

$$\|\hat{\lambda} - \lambda\|_{\ell_2} \leq 2\left(\sum_{j \in J_0 \cup J_1} |\hat{\lambda}_j - \lambda_j|^2\right)^{1/2} \leq 2\beta_2(d(\lambda))\|f_{\hat{\lambda}} - f_\lambda\|_{L_2(\Pi)}. \quad (2.9)$$

We will also use Lemma 6 to bound

$$\sup_{\|u\|_{\ell_2} \leq 1, u \in C_{J_\lambda}} |(\Pi_n - \Pi)(|f_u|)| \leq C \sup_{\|u\|_{\ell_2} \leq 1, d(u) \leq \bar{d}} \|f_u\|_{\psi_1} \sqrt{\frac{A\bar{d}\log(N/\bar{d})}{n}}, \quad (2.10)$$

which holds with probability at least

$$1 - 5^{-\bar{d}A}\binom{N}{\leq \bar{d}}^{-A}.$$

The first term in the right-hand side of (2.8) is bounded as in the proof of Theorem 1

$$\|f_{\hat{\lambda}} - f_\lambda\|_{L_1(\Pi_n)} \leq \sqrt{2\varepsilon\|\hat{\lambda} - \lambda\|_{\ell_1}} \quad (2.11)$$

and we then use

$$\begin{aligned}\|\hat{\lambda} - \lambda\|_{\ell_1} &\leq 2\sum_{j \in J}|\hat{\lambda}_j - \lambda_j| \leq 2\sqrt{d(\lambda)}\left(\sum_{j \in J \cup J_1}|\hat{\lambda}_j - \lambda_j|^2\right)^{1/2} \\ &\leq 2\beta_2(d(\lambda))\sqrt{d(\lambda)}\|f_{\hat{\lambda}} - f_\lambda\|_{L_2(\Pi)}.\end{aligned} \quad (2.12)$$

It remains to substitute bounds (2.9)–(2.12) into (2.8), to also use (2.2) and to solve the resulting inequality with respect to $\|f_{\hat{\lambda}} - f_\lambda\|_{L_2(\Pi)}$ to obtain the first bound of the theorem.

To prove the second bound, it is enough to use (2.9) to obtain a bound on $\|\hat{\lambda} - \lambda\|_{\ell_2}$. The remaining two bounds are proved the same way as in the proof of Theorem 1. □

Condition (2.2) required in Theorem 2 is rather restrictive. Moreover, since the $\psi_1$-norm of $f_u$ is involved in the definition of the set $\Lambda_{S,2}(A)$ of sparse vectors, one might need the equivalence of the $\psi_1$- and the $L_2$-norms on the linear span of the dictionary in order to have a more explicit way to describe the sparsity of the problem. Condition (2.2) is not needed in Theorem 1. However, this condition is needed to bound the quantity $\beta(J)$ in terms of the quantity $\beta_2(J)$, the latter being much more convenient because of its simple relationships to various geometric characteristics of the dictionary (see Section 3). So, in both cases, one must rely on condition (2.2) and the class of examples to which the results apply is rather limited (such as Gaussian and Rademacher dictionaries). Below, we give another version of sparsity oracle inequalities for the Dantzig selector that does not have this drawback and which applies to a variety of dictionaries. However, in this case, much more is required in terms of sparsity. In the orthonormal case, the result applies only to oracle vectors $\lambda$ with

$$d(\lambda) \leq c\sqrt{\frac{n}{\log N}}.$$

A similar constraint was needed, for instance, in sparsity oracle inequalities for LASSO in the paper by Bunea, Tsybakov and Wegkamp (2007). In Theorems 1 and 2, the oracle sets were larger, including vectors $\lambda$ with $d(\lambda)$ comparable to $n$.

The set of "sparse" vectors is defined as

$$\Lambda_{S,3}(A) := \left\{ \lambda \in \mathbb{R}^N : C\beta_2^2(J_\lambda)d(\lambda) \max_{1 \leq k,j \leq N} \|h_k h_j\|_{\psi_1} \sqrt{\frac{A \log N}{n}} \leq 1/8 \right\}$$

and the oracle set becomes

$$\tilde{\Lambda}_\varepsilon^3(A) := \Lambda_\varepsilon(A) \cap \Lambda_{S,3}(A).$$

**Theorem 3.** *There exists a constant $C$ in the definitions of $\Lambda_\varepsilon(A)$ and $\Lambda_{S,3}(A)$ such that for all $A \geq 1$ with probability at least $1 - N^{-A}$, the following bounds hold: $\forall \lambda \in \hat{\Lambda}_\varepsilon \cap \Lambda_{S,3}(A)$*

$$\|f_{\hat{\lambda}} - f_\lambda\|_{L_2(\Pi)} \leq 8\beta_2(J_\lambda)\sqrt{d(\lambda)}\varepsilon,$$

$$\|\hat{\lambda} - \lambda\|_{\ell_1} \leq 16\beta_2^2(J_\lambda)d(\lambda)\varepsilon$$

*and*

$$\|\hat{\lambda} - \lambda\|_{\ell_2} \leq 16\beta_2^2(d(\lambda))\sqrt{d(\lambda)}\varepsilon.$$

*Suppose that $\tilde{\Lambda}_\varepsilon^3(A) \neq \varnothing$. Then, with probability at least $1 - N^{-A}$,*

$$\|f_{\hat{\lambda}} - f_*\|_{L_2(\Pi)} \leq \inf_{\lambda \in \tilde{\Lambda}_\varepsilon^3(A)} [\|f_\lambda - f_*\|_{L_2(\Pi)} + 8\beta_2(J_\lambda)\sqrt{d(\lambda)}\varepsilon].$$

*Moreover, if, in addition, $f_* = f_{\lambda^*}, \lambda^* \in \mathbb{R}^N$, then, also with the same probability,*

$$\|\hat{\lambda} - \lambda^*\|_{\ell_1} \leq \inf_{\lambda \in \tilde{\Lambda}_\varepsilon^3(A)} [\|\lambda - \lambda^*\|_{\ell_1} + 16\beta_2^2(J_\lambda)d(\lambda)\varepsilon]$$

*and*

$$\|\hat{\lambda} - \lambda^*\|_{\ell_2} \leq \inf_{\lambda \in \tilde{\Lambda}_\varepsilon^3(A)} [\|\lambda - \lambda^*\|_{\ell_2} + 16\beta_2^2(d(\lambda))\sqrt{d(\lambda)}\varepsilon].$$

**Proof.** This is similar to the proofs of Theorems 1 and 2. The following bounds are used for all $\lambda \in \hat{\Lambda}_\varepsilon \cap \Lambda_{S,3}(A)$:

$$\|\hat{\lambda} - \lambda\|_{\ell_1} \leq 2\beta_2(J_\lambda)\sqrt{d(\lambda)}\|f_{\hat{\lambda}} - f_\lambda\|_{L_2(\Pi)},$$
$$\|\hat{\lambda} - \lambda\|_{\ell_2} \leq 2\beta_2(d(\lambda))\|f_{\hat{\lambda}} - f_\lambda\|_{L_2(\Pi)},$$

and

$$\|f_{\hat{\lambda}} - f_\lambda\|^2_{L_2(\Pi)} \leq \|f_{\hat{\lambda}} - f_\lambda\|^2_{L_2(\Pi_n)} + (\Pi - \Pi_n)(|f_{\hat{\lambda}} - f_\lambda|^2)$$
$$\leq 2\varepsilon\|\hat{\lambda} - \lambda\|_{\ell_1} + \sup_{\|u\|_{\ell_1} \leq 1}|(\Pi_n - \Pi)(|f_u|^2)|\|\hat{\lambda} - \lambda\|^2_{\ell_1},$$

and then Lemma 5 is applied to bound the last term on the right-hand side. $\square$

It is not our goal in this paper to study the fixed design case in detail. However, some results are rather easy to obtain, in a manner similar to our derivations in the random design case (actually, with some simplifications). In particular, the following result holds. We will use a version of $\beta_2(J)$ with $\Pi$ replaced by the empirical measure $\Pi_n$ (based on the design points):

$$\hat{\beta}_2(J) := \beta_2(J; \Pi_n) \quad \text{and} \quad \hat{\beta}_2(d) := \beta_2(d; \Pi_n).$$

**Theorem 4.** *Suppose $X_1, \ldots, X_n$ are non-random design points in $S$ and let $\Pi_n$ be the empirical measure based on $X_1, \ldots, X_n$. Suppose, also, that $f_* = f_{\lambda^*}$, $\lambda^* \in \mathbb{R}^N$. There exists a constant $C > 0$ such that, for all $A \geq 1$ and all*

$$\varepsilon \geq C\|\xi\|_{\psi_2} \max_{1 \leq k \leq N} \|h_k\|_{L_2(\Pi_n)}\sqrt{\frac{A\log N}{n}},$$

*with probability at least $1 - N^{-A}$, the following bounds hold:*

$$\|f_{\hat{\lambda}} - f_{\lambda^*}\|_{L_2(\Pi_n)} \leq 4\hat{\beta}_2(J_{\lambda^*})\sqrt{d(\lambda^*)}\varepsilon,$$
$$\|\hat{\lambda} - \lambda^*\|_{\ell_1} \leq 8\hat{\beta}_2^2(J_{\lambda^*})d(\lambda^*)\varepsilon$$

*and*

$$\|\hat{\lambda} - \lambda^*\|_{\ell_2} \leq 8\hat{\beta}_2^2(d(\lambda^*))\sqrt{d(\lambda^*)}\varepsilon.$$

**Proof.** This is, essentially, a simplified version of the arguments used in the proofs of Theorems 1 and 2. The following two bounds are obtained exactly as in the proof of Theorem 1:

$$\|f_{\hat\lambda} - f_{\lambda^*}\|_{L_2(\Pi_n)} \leq \sqrt{2\varepsilon \|\hat\lambda - \lambda^*\|_{\ell_1}} \tag{2.13}$$

and

$$\|\hat\lambda - \lambda^*\|_{\ell_1} \leq 2\hat\beta_2(J_{\lambda^*})\sqrt{d(\lambda^*)}\|f_{\hat\lambda} - f_{\lambda^*}\|_{L_2(\Pi_n)}. \tag{2.14}$$

They hold if $\lambda^* \in \hat\Lambda_\varepsilon$, which is equivalent to the condition

$$\max_{1 \leq k \leq N} \left| n^{-1} \sum_{j=1}^n \xi_j h_k(X_j) \right| \leq \varepsilon.$$

If $\|\xi\|_{\psi_2} < +\infty$ and

$$\varepsilon \geq C\|\xi\|_{\psi_2} \max_{1 \leq k \leq N} \|h_k\|_{L_2(\Pi_n)} \sqrt{\frac{A \log N}{n}},$$

then standard exponential bounds for sums of independent random variables and bounds on the maximum of random variables in Orlicz spaces imply that with probability at least $1 - N^{-A}$, $\lambda^* \in \hat\Lambda_\varepsilon$.

It remains to combine (2.13) and (2.14) to prove that with probability at least $1 - N^{-A}$,

$$\|f_{\hat\lambda} - f_{\lambda^*}\|_{L_2(\Pi_n)} \leq 4\hat\beta_2(J_{\lambda^*})\sqrt{d(\lambda^*)}\varepsilon$$

and

$$\|\hat\lambda - \lambda^*\|_{\ell_1} \leq 8\hat\beta_2^2(J_{\lambda^*})d(\lambda^*)\varepsilon.$$

Arguing exactly as in the proof of Theorem 2 (in particular, using Lemma 1), one can add to this that

$$\|\hat\lambda - \lambda^*\|_{\ell_2} \leq 8\hat\beta_2^2(d(\lambda^*))\sqrt{d(\lambda^*)}\varepsilon. \qquad \square$$

One can also obtain an upper bound on $\hat\beta_2(J), d(J) = d$ in terms of fixed design versions of restricted isometry constants (see Lemma 2, Corollary 6), which leads to Theorem 1 in Candes and Tao (2007) (bounds on the performance of the Dantzig selector under the UUP). They also proved a sharper oracle inequality in the case of random design that we are not going to reproduce here. However, Corollary 3 in the next section provides a direct proof of a similar inequality in the random design case (for orthonormal dictionaries).

## 3. Corollaries and remarks

Under the additional assumption

$$\|\lambda\|_{\ell_2} \leq B\|f_\lambda\|_{L_1(\Pi)}, \qquad \lambda \in \mathbb{R}^N, \tag{3.1}$$

it is easy to establish a corollary of Theorem 1 that implies the main results of Candes and Tao (2007) in the random design case. Assume that $f_* = f_{\lambda^*}$, $\lambda^* \in \mathbb{R}^N$.

Define

$$\bar{\Lambda}_\varepsilon(A) := \Lambda_\varepsilon(A) \cap \left\{ \lambda \in \mathbb{R}^N : CB \max_{1 \leq k \leq N} \|h_k(X)\|_{\psi_1} \sqrt{\frac{Ad(\lambda)\log N}{n}} \leq 1/4 \right\}.$$

**Corollary 1.** *Suppose that condition (3.1) holds. There then exists a constant $C$ in the definition of the set $\bar{\Lambda}_\varepsilon(A)$ such that, for all $A \geq 1$ and under the assumption $\bar{\Lambda}_\varepsilon(A) \neq \varnothing$, with probability at least $1 - N^{-A}$,*

$$\|\hat{\lambda} - \lambda^*\|_{\ell_2} \leq \inf_{\lambda \in \bar{\Lambda}_\varepsilon(A)} [\|\lambda - \lambda^*\|_{\ell_2} + 16B^2 \sqrt{d(\lambda)}\varepsilon].$$

**Proof.** Under the assumption (3.1),

$$\sum_{j \in J} |\lambda_j| \leq \sqrt{d(J)} \|\lambda\|_{\ell_2} \leq B\sqrt{d(J)} \|f_\lambda\|_{L_1(\Pi)},$$

implying that $\beta(J) \leq B\sqrt{d(J)}$. Therefore, $\bar{\Lambda}_\varepsilon(A) \subset \tilde{\Lambda}_\varepsilon(A)$. Denoting by $\bar{\lambda}$ the value of $\lambda$ that minimizes the right-hand side of the bound of Corollary 1, the bound

$$\|f_{\hat{\lambda}} - f_{\bar{\lambda}}\|_{L_1(\Pi)} \leq 16\beta(J_{\bar{\lambda}})\varepsilon \leq 16B\sqrt{d(\bar{\lambda})}\varepsilon$$

follows from the first inequality of Theorem 1. This yields

$$\|\hat{\lambda} - \bar{\lambda}\|_{\ell_2} \leq 16B^2 \sqrt{d(\bar{\lambda})}\varepsilon,$$

implying the result. □

If

$$\varepsilon \geq C \max_{1 \leq k \leq N} \|\xi h_k(X)\|_{\psi_1} \sqrt{\frac{A \log N}{n}} \tag{3.2}$$

and the vector $\lambda^*$ is sufficiently sparse in the sense that

$$CB \max_{1 \leq k \leq N} \|h_k(X)\|_{\psi_1} \sqrt{\frac{Ad(\lambda^*)\log N}{n}} \leq 1/4, \tag{3.3}$$

then Corollary 1 immediately implies that with probability at least $1 - N^{-A}$,

$$\|\hat{\lambda} - \lambda^*\|_{\ell_2} \leq 16B^2 \sqrt{d(\lambda^*)}\varepsilon.$$

It is enough to observe that, in this case, $\lambda^* \in \bar{\Lambda}_\varepsilon(A)$ and to use $\lambda = \lambda^*$ in the bound of the corollary (without taking the infimum). By simple properties of Orlicz norms,

$$\|\xi h_k(X)\|_{\psi_1} \leq \|\xi\|_{\psi_2} \|h_k(X)\|_{\psi_2}. \tag{3.4}$$

(Indeed, for random variables $\eta_1, \eta_2$ such that $\|\eta_i\|_{\psi_2} \leq 1$, the following holds by the definitions of the norms:

$$\|\eta_1\eta_2\|_{\psi_1} \leq \|(\eta_1^2 + \eta_2^2)/2\|_{\psi_1} \leq (\|\eta_1^2\|_{\psi_1} + \|\eta_2^2\|_{\psi_1})/2 \leq (\|\eta_1\|_{\psi_2} + \|\eta_2\|_{\psi_2})/2 \leq 1.$$

This immediately implies that for all $\eta_1, \eta_2$,

$$\|\eta_1\eta_2\|_{\psi_1} \leq \|\eta_1\|_{\psi_2}\|\eta_2\|_{\psi_2}.)$$

If $\xi$ is a normal random variable with mean zero and variance $\sigma^2$, the $\psi_2$-norm of $\xi$ coincides with $\sigma$ (up to a numerical constant). So, under the assumption that

$$\|h_k(X)\|_{\psi_2} \leq 1, \qquad k = 1, \ldots, N,$$

conditions (3.2) and (3.3) take the following form:

$$\varepsilon \geq C\sigma\sqrt{\frac{A\log N}{n}} \tag{3.5}$$

and

$$CB\sqrt{\frac{Ad(\lambda^*)\log N}{n}} \leq 1/4. \tag{3.6}$$

The case $\sigma = 0$ (no noise in the regression model) is of special interest. In this case, $\|\xi h_k(X)\|_{\psi_1} = 0$ and one can use $\varepsilon = 0$ in the definition of the Dantzig selector. The following result holds.

**Corollary 2.** *Suppose that $\xi = 0$ and $\|h_k(X)\|_{\psi_1} \leq 1$. Let $\varepsilon = 0$. If condition (3.1) and sparsity condition (3.6) hold, then, with probability at least $1 - N^{-A}$, $\hat{\lambda} = \lambda^*$.*

Moreover, if we assume that both $\psi_1$- and $L_1$-norms on the linear span of the dictionary are equivalent, up to numerical constants (independent of $N$), to the $\ell_2$-norm in the space of vectors of coefficients (which is true, for instance, for Gaussian and Rademacher dictionaries), then the sparsity assumption (3.6) can be replaced by a slightly weaker assumption $d(\lambda^*) \leq \bar{d}$, where $\bar{d}$ satisfies the condition

$$\sqrt{\frac{A\bar{d}\log(N/\bar{d})}{n}} \leq c, \tag{3.7}$$

with a proper choice of $c$, and the conclusion of Corollary 2 still holds (when $N = n$, this means that $d(\lambda^*) \leq cn$, with a proper choice of constant $c$). Theorem 2 must be used, leading to the bound

$$\mathbb{P}\{\hat{\lambda} \neq \lambda^*\} \leq 5^{-\bar{d}A}\binom{N}{\leq \bar{d}}^{-A},$$

so, in this case, the probability of the error is bounded in a better way.

Hence, the Dantzig selector provides an *exact* solution to the problem of recovery of a sparse vector $\lambda^*$ based on noiseless measurements of function $f_{\lambda^*}$ at random points. It is easy to see that in this case, one can also use another definition of $\hat{\lambda}$, as a minimizer of the $\ell_1$-norm $\|\lambda\|_{\ell_1}$ subject to the linear constraints $f_\lambda(X_j) = Y_j$, $j = 1, \ldots, N$ with no changes in the proof. This striking fact has been known for a while and has some interesting connections to deep results in convex geometry and asymptotic geometric analysis (such as, for instance, neighborliness of convex polytopes; see Donoho (2006a, 2006b), Candes and Tao (2005), Candes, Romberg and Tao (2006), Rudelson and Vershynin (2005), Mendelson, Pajor and Tomczak-Jaegermann (2007) and references therein).

We will consider another interesting consequence of Corollary 1 under the additional assumptions that the dictionary is orthonormal in $L_2(\Pi)$ and that the $\psi_2$- and the $L_2$-norms are equivalent on the linear span of the dictionary up to a numerical constant. Because of orthonormality, the $L_2$-norm is equal to the $\ell_2$-norm in the space of coefficients, and condition (3.1) becomes, in this case, a version of condition (2.2):

$$\|f_\lambda\|_{L_2(\Pi)} \leq B \|f_\lambda\|_{L_1(\Pi)}, \qquad \lambda \in \mathbb{R}^N.$$

Thus, all of the Orlicz norms between $L_1$ and $\psi_2$ are equivalent in this case. In particular, this applies to Gaussian and Rademacher dictionaries. The result given below also follows from the oracle inequality proven by Candes and Tao (2007) in the fixed design case (under the UUP condition). It is also not hard to establish it for the dictionaries that are not necessarily orthonormal, but that satisfy some assumption on the "weakness" of correlations between functions $h_j$.

**Corollary 3.** *Under the above assumptions, including (3.1), the assumption that the noise is $N(0; \sigma^2)$, that the dictionary is orthonormal in $L_2(\Pi)$ and that the $\psi_2$- and the $L_2$-norms are equivalent on the linear span of the dictionary up to a numerical constant, there exists a choice of constant $C$ such that for all $\varepsilon$ satisfying condition (3.5) and $\lambda^*$ satisfying the sparsity assumption (3.6), with probability at least $1 - N^{-A}$ and with some $D > 0$ depending on $B$ in condition (3.1),*

$$\|\hat{\lambda} - \lambda^*\|_{\ell_2}^2 \leq D \sum_{j=1}^{N} (|\lambda_j^*|^2 \wedge \varepsilon^2) = D \inf_{J \subset \{1, \ldots, N\}} \left[ \sum_{j \notin J} |\lambda_j^*|^2 + d(J)\varepsilon^2 \right].$$

**Proof.** Define $\bar{\lambda}^*$ as follows:

$$\bar{\lambda}_j^* = \lambda_j^* I(|\lambda_j^*| \geq \varepsilon/3).$$

We then have

$$|\langle f_{\bar{\lambda}^*} - f_{\lambda^*}, h_k \rangle_{L_2(\Pi)}| = |\lambda_k^*| \leq \varepsilon/3$$

for all $k \in J_{\lambda^*}$, $|\lambda_k^*| \leq \varepsilon/3$ and it is equal to 0 otherwise.

We also have, for all $k = 1, \ldots, N$,

$$\|\xi h_k(X)\|_{\psi_1} \leq \|\xi\|_{\psi_2} \|h_k(X)\|_{\psi_2} \leq c\sigma,$$

with a numerical constant $c$, implying that with a proper choice of $C, C'$,

$$C\|\xi h_k(X)\|_{\psi_1}\sqrt{\frac{A\log N}{n}} \leq \varepsilon/3,$$

provided that

$$\varepsilon \geq C'\sigma\sqrt{\frac{A\log N}{n}}.$$

Finally,

$$\|(f_{\bar{\lambda}^*} - f_{\lambda^*})(X)h_k(X)\|_{\psi_1} \leq \|(f_{\bar{\lambda}^*} - f_{\lambda^*})(X)\|_{\psi_2}\|h_k(X)\|_{\psi_2} \leq c\|f_{\bar{\lambda}^*} - f_{\lambda^*}\|_{L_2(\Pi)}$$

$$= c\bigg(\sum_{j\,:\,|\lambda_j^*|<\varepsilon/3}|\lambda_j^*|^2\bigg)^{1/2} \leq c(\varepsilon/3)\sqrt{d(\lambda^*)},$$

which implies that

$$C\|(f_{\bar{\lambda}^*} - f_{\lambda^*})(X)h_k(X)\|_{\psi_1}\sqrt{\frac{A\log N}{n}} \leq \varepsilon/3,$$

provided that

$$cC\sqrt{\frac{Ad(\lambda^*)\log N}{n}} \leq 1.$$

The last condition is equivalent to (3.6) with a proper choice of constant therein.

Hence, $\bar{\lambda}^* \in \bar{\Lambda}_\varepsilon(A)$ and Corollary 1 implies that with probability at least $1 - N^{-A}$,

$$\|\hat{\lambda} - \lambda^*\|_{\ell_2} \leq \bigg(\sum_{j:|\lambda_j^*|<\varepsilon/3}|\lambda_j^*|^2\bigg)^{1/2} + 16B^2\sqrt{\mathrm{card}(j:|\lambda_j^*|\geq \varepsilon/3)}\varepsilon,$$

which yields that, with some constant $D$ (depending on $B$),

$$\|\hat{\lambda} - \lambda^*\|_{\ell_2}^2 \leq D\sum_{j=1}^{N}(|\lambda_j^*|^2 \wedge \varepsilon^2). \qquad \square$$

We now describe a couple of ways of bounding the quantity $\beta_2(J)$ involved in Theorem 2 and in the upper bound on $\beta(J)$.

Let $\kappa(J)$ denote the minimal eigenvalue of the Gram matrix $(\langle h_i, h_j\rangle_{L_2(\Pi)})_{i,j\in J}$. Also, denote by $L_J$ the linear span of $\{h_j : j \in J\}$ and let

$$\rho(J) := \sup_{f\in L_J, g\in L_{J^c}, f,g\neq 0}\bigg|\frac{\langle f,g\rangle_{L_2(\Pi)}}{\|f\|_{L_2(\Pi)}\|g\|_{L_2(\Pi)}}\bigg|.$$

$\rho(J)$ is thus the largest "correlation coefficient" (or the largest cosine of the angle) between functions in the linear span of a subset $\{h_j : j \in J\}$ of the dictionary and the linear span of its complement. It is of interest to compare $\rho(J)$ with the concept of canonical correlation used in multivariate statistical analysis. It is very easy to check (see Koltchinskii (2009), Proposition 1) that

$$\beta_2(J) \leq \frac{1}{\sqrt{\kappa(J)(1-\rho^2(J))}}$$

and, as a consequence, if (2.2) holds with some constant $B(J) > 0$, then

$$\beta(J) \leq \sqrt{\tilde{d}(J)},$$

where

$$\tilde{d}(J) := \frac{B^2(J)d(J)}{\kappa(J)(1-\rho^2(J))}.$$

In particular, if the dictionary $\{h_1, \ldots, h_N\}$ is orthonormal in $L_2(\Pi)$ and condition (2.2) holds with a constant $B$ that does not depend on $J$ (for instance, in the case of Gaussian or Bernoulli dictionaries), then $\kappa(J) = 1$ and $\rho(J) = 0$, so $\tilde{d}(J) = B^2 d(J)$, leading to the bound

$$\beta(J) \leq B\sqrt{d(J)}.$$

We define

$$\tilde{d}(\lambda) := \tilde{d}(J_\lambda),$$

which plays the role of a modified "dimension" of the vector $\lambda$ (that takes into account how close the dictionary is to the orthonormality property; in the orthonormal case, $\tilde{d}(\lambda) = B^2 d(\lambda)$).

Define

$$\check{\Lambda}_\varepsilon(A) := \Lambda_\varepsilon(A) \cap \left\{\lambda \in \mathbb{R}^N : C \max_{1 \leq k \leq N} \|h_k(X)\|_{\psi_1} \sqrt{\frac{A\tilde{d}(\lambda)\log N}{n}} \leq 1/4\right\}.$$

The proof of the following corollary repeats the proof of Corollary 1, with the $\ell_2$-norm replaced by the $L_2(\Pi)$-norm.

**Corollary 4.** *Suppose that condition (2.2) holds. There exists a constant $C$ in the definition of $\check{\Lambda}_\varepsilon(A)$ such that for all $A \geq 1$, the assumption $\check{\Lambda}_\varepsilon(A) \neq \varnothing$ implies that with probability at least $1 - N^{-A}$,*

$$\|f_{\hat{\lambda}} - f_*\|_{L_2(\Pi)} \leq \inf_{\lambda \in \check{\Lambda}_\varepsilon(A)} [\|f_\lambda - f_*\|_{L_2(\Pi)} + 16\sqrt{\tilde{d}(\lambda)}\varepsilon].$$

Denote by $f_{\lambda^0}$, $\lambda^0 \in \mathbb{R}^N$ the orthogonal projection in $L_2(\Pi)$ of $f_*$ onto the linear span of the dictionary. The following result is also straightforward.

**Corollary 5.** *Suppose that the condition (2.2) holds and that the noise is normal with mean zero and variance $\sigma^2$. There then exists a constant $C$ such that for all $A \geq 1$ and all*

$$\varepsilon \geq C(\|f_{\lambda^0} - f_*\|_{\psi_2} + \sigma)\sqrt{\frac{A \log N}{n}},$$

*if $\lambda^0$ satisfies the "sparsity" condition*

$$C\sqrt{\frac{A\tilde{d}(\lambda^0) \log N}{n}} \leq 1/4,$$

*then with probability at least $1 - N^{-A}$,*

$$\|f_{\hat{\lambda}} - f_*\|^2_{L_2(\Pi)} \leq \|f_{\lambda^0} - f_*\|^2_{L_2(\Pi)} + 16^2 \tilde{d}(\lambda^0)\varepsilon^2.$$

**Proof.** Under the assumptions, it is easy to check that $\lambda^0 \in \check{\Lambda}_\varepsilon(A)$. This allows one to deduce that, with probability at least $1 - N^{-A}$,

$$\|f_{\hat{\lambda}} - f_{\lambda^0}\|^2_{L_2(\Pi)} \leq 16^2 \tilde{d}(\lambda^0)\varepsilon^2.$$

Since $f_{\hat{\lambda}} - f_{\lambda^0}$ and $f_{\lambda^0} - f_*$ are orthogonal, this implies the result. □

Another approach to bounding $\beta_2(J)$ is based on some quantities involved in the restricted isometry condition.

For $I, J \subset \{1, \ldots, N\}$, $I \cap J = \varnothing$, define

$$r(I; J) := \sup_{f \in L_I, g \in L_J, f, g \neq 0} \left| \frac{\langle f, g \rangle_{L_2(\Pi)}}{\|f\|_{L_2(\Pi)} \|g\|_{L_2(\Pi)}} \right|$$

(for $\rho(J)$ defined before, $\rho(J) = r(J, J^c)$). Let

$$\rho_d := \max\{r(I, J) : I, J \subset \{1, \ldots, N\},\ I \cap J = \varnothing,\ d(I) = 2d,\ d(J) = d\}.$$

This quantity measures the correlation between linear spans of disjoint parts of the dictionary of fixed cardinalities, $d$ and $2d$ (it is a more "local" characteristic of the dictionary than the quantity $\rho(J)$ used before).

We will also define

$$m_d := \inf\{\|f_u\|_{L_2(\Pi)} : u \in \mathbb{R}^N, \|u\|_{\ell_2} = 1, d(u) \leq d\}$$

and

$$M_d := \sup\{\|f_u\|_{L_2(\Pi)} : u \in \mathbb{R}^N, \|u\|_{\ell_2} = 1, d(u) \leq d\}.$$

**Lemma 2.** *Suppose $J \subset \{1, \ldots, N\}$, $d(J) = d$ and*

$$\rho_d < \frac{m_{2d}}{M_{2d}}.$$

*Then,*

$$\beta_2(J) \leq \frac{1}{m_{2d} - \rho_d M_{2d}}.$$

**Proof.** Recall Lemma 1 and its notation. Denote by $P_I$ the orthogonal projection on $L_I \subset L_2(\Pi)$. Then, for all $u \in C_J$,

$$\left\| \sum_{j=1}^{N} u_j h_j \right\|_{L_2(\Pi)} \geq \left\| P_{J_0 \cup J_1} \sum_{j=1}^{N} u_j h_j \right\|_{L_2(\Pi)}$$

$$\geq \left\| \sum_{j \in J_0 \cup J_1} u_j h_j \right\|_{L_2(\Pi)} - \left\| P_{J_0 \cup J_1} \sum_{j \notin J_0 \cup J_1} u_j h_j \right\|_{L_2(\Pi)}$$

$$\geq \left\| \sum_{j \in J_0 \cup J_1} u_j h_j \right\|_{L_2(\Pi)} - \sum_{k \geq 2} \left\| P_{J_0 \cup J_1} \sum_{j \in J_k} u_j h_j \right\|_{L_2(\Pi)}$$

$$\geq \left\| \sum_{j \in J_0 \cup J_1} u_j h_j \right\|_{L_2(\Pi)} - \rho_d \sum_{k \geq 2} \left\| \sum_{j \in J_k} u_j h_j \right\|_{L_2(\Pi)}$$

$$\geq \left\| \sum_{j \in J_0 \cup J_1} u_j h_j \right\|_{L_2(\Pi)} - \rho_d M_{2d} \sum_{k \geq 2} \| u^{(k)} \|_{\ell_2}$$

$$\geq \left\| \sum_{j \in J_0 \cup J_1} u_j h_j \right\|_{L_2(\Pi)} - \rho_d M_{2d} \left( \sum_{j \in J \cup J_1} |u_j|^2 \right)^{1/2}$$

$$\geq \left\| \sum_{j \in J_0 \cup J_1} u_j h_j \right\|_{L_2(\Pi)} - \rho_d \frac{M_{2d}}{m_{2d}} \left\| \sum_{j \in J_0 \cup J_1} u_j h_j \right\|_{L_2(\Pi)}$$

$$= \left( 1 - \rho_d \frac{M_{2d}}{m_{2d}} \right) \left\| \sum_{j \in J_0 \cup J_1} u_j h_j \right\|_{L_2(\Pi)}.$$

On the other hand, for $u \in C_J$,

$$\left( \sum_{j \in J} |u_j|^2 \right)^{1/2} \leq \left( \sum_{j \in J_0 \cup J_1} |u_j|^2 \right)^{1/2} \leq m_{2d}^{-1} \left\| \sum_{j \in J_0 \cup J_1} u_j h_j \right\|_{L_2(\Pi)},$$

implying that

$$\left(\sum_{j\in J}|u_j|^2\right)^{1/2} \leq m_{2d}^{-1}\left(1-\rho_d\frac{M_{2d}}{m_{2d}}\right)^{-1}\left\|\sum_{j=1}^N u_j h_j\right\|_{L_2(\Pi)},$$

which yields

$$\beta_2(J) \leq \frac{1}{m_{2d}-\rho_d M_{2d}}. \qquad \square$$

If $m_d \leq 1 \leq M_d \leq 2$, one can express the restricted isometry constant $\delta_d = \delta_d(\Pi)$ as

$$\delta_d = (M_d - 1) \vee (1 - m_d).$$

It is also easy to show that

$$\rho_d \leq \left[\left(\frac{1+\delta_{3d}}{1-\delta_{3d}}\right)^2 - 1\right] \vee \left[1 - \left(\frac{1-\delta_{3d}}{1+\delta_{3d}}\right)^2\right].$$

Lemma 2 then implies that there exists $\bar{\delta} < 1$ such that, under the condition $\delta_{3d} \leq \bar{\delta}$, $\beta_2(J) \leq \bar{c}$ for all sets $J$ with $d(J) \leq d$, where $\bar{c}$ is a constant that depends only on $\bar{\delta}$ (for instance, one can take $\bar{\delta} = 1/8$).

Denote by $\bar{d}$ the largest $d$ such that $d \leq N/e - 1, \delta_{3d} \leq \bar{\delta}$,

$$\frac{A d \log(N/d)}{n} \leq 1$$

and

$$CB(d)\sup_{\|u\|_{\ell_2}\leq 1, d(u)\leq d}\|f_u\|_{\psi_1}\sqrt{\frac{Ad\log(N/d)}{n}} \leq 1/4.$$

Let

$$\check{\Lambda}_\varepsilon^2(A) := \Lambda_\varepsilon(A) \cap \{\lambda \in \mathbb{R}^N : d(\lambda) \leq \bar{d}\}.$$

The next corollary is an immediate consequence of Theorem 2 and Lemma 2. It shows that sparse enough target functions can be recovered by the Dantzig selector under a version of the restricted isometry assumption. In particular, Proposition 1 in the Introduction immediately follows from this corollary.

**Corollary 6.** *There exist constants $C, D$ depending only on $\bar{\delta}$ such that for all $A \geq 1$, the assumption $\check{\Lambda}_\varepsilon^2(A) \neq \varnothing$ implies that with probability at least $1 - N^{-A}$,*

$$\|f_{\hat{\lambda}} - f_*\|_{L_2(\Pi)} \leq \inf_{\lambda \in \check{\Lambda}_\varepsilon^2(A)}[\|f_\lambda - f_*\|_{L_2(\Pi)} + DB^2(d(\lambda))\sqrt{d(\lambda)}\varepsilon].$$

Moreover, if, in addition, $f_* = f_{\lambda^*}, \lambda^* \in \mathbb{R}^N$, then we also have

$$\|\hat{\lambda} - \lambda^*\|_{\ell_2} \leq \inf_{\lambda \in \tilde{\Lambda}^2_\varepsilon(A)}[\|\lambda - \lambda^*\|_{\ell_2} + DB^2(d(\lambda))\sqrt{d(\lambda)}\varepsilon].$$

Another way to bound the quantity $\beta_2(J)$ is given in the following proposition that is akin to some statements in Bickel, Ritov and Tsybakov (2009) (in the fixed design case). The proof is rather straightforward and is based on a simple modification of Lemma 1.

**Proposition 3.** *If $J \subset \{1, \ldots, N\}$ with $d(J) = d$ and, for some $s \geq 1$,*

$$\frac{M_s}{m_{s+d}} < \sqrt{\frac{s}{d}},$$

*then*

$$\beta_2(J) \leq \frac{\sqrt{s}}{\sqrt{s}m_{d+s} - \sqrt{d}M_s}.$$

We conclude this section with a couple of examples that provide some explanation of the role of such geometric characteristics of the dictionary as $\beta_2(J)$ in sparse recovery problems.

***Example.*** Consider the case where the functions $h_1, \ldots, h_N$ are orthogonal in $L_2(\Pi)$. It is easy to see that

$$\beta_2(J) = \frac{1}{\min_{j \in J} \|h_j\|_{L_2(\Pi)}}.$$

Suppose that $f_* = f_{\lambda^*}$ with $\lambda^* \in \mathbb{R}^N$ and

$$\|h_j\|_{L_2(\Pi)} = \tau > 0, \qquad j \in J_{\lambda^*}.$$

Fix the value of the parameter $\varepsilon > 0$ of the Dantzig selector and consider, for simplicity, the limit case when $n \to \infty$. In this limit, the set $\hat{\Lambda}_\varepsilon$ becomes

$$\hat{\Lambda}_\varepsilon := \{\lambda \in \mathbb{R}^N : |\langle f_\lambda - f_{\lambda^*}, h_k \rangle| \leq \varepsilon, k = 1, \ldots, N\}$$

which, under the orthogonality assumption, is just

$$\hat{\Lambda}^\varepsilon := \{\lambda \in \mathbb{R}^N : |\lambda_k - \lambda^*_k|\|h_k\|^2_{L_2(\Pi)} \leq \varepsilon, k = 1, \ldots, N\}.$$

It is easy to compute the Dantzig selector: for $k = 1, \ldots, N$,

$$\hat{\lambda}_k = \left(\lambda^*_k - \frac{\varepsilon}{\|h_k\|^2_{L_2(\Pi)}}\right)I\left(\lambda^*_k \geq \frac{\varepsilon}{\|h_k\|^2_{L_2(\Pi)}}\right) + \left(\lambda^*_k + \frac{\varepsilon}{\|h_k\|^2_{L_2(\Pi)}}\right)I\left(\lambda^*_k \leq -\frac{\varepsilon}{\|h_k\|^2_{L_2(\Pi)}}\right).$$

Therefore,

$$\|\hat{\lambda} - \lambda^*\|_{\ell_2}^2 = \varepsilon^2 \sum_{k:|\lambda_k^*| \geq \varepsilon/\|h_k\|_{L_2(\Pi)}^2} \frac{1}{\|h_k\|_{L_2(\Pi)}^4} + \sum_{k:|\lambda_k^*| < \varepsilon/\|h_k\|_{L_2(\Pi)}^2} |\lambda_k^*|^2.$$

If $|\lambda_j^*| \geq \varepsilon/\tau^2$ for all $j \in J_{\lambda^*}$, we get

$$\|\hat{\lambda} - \lambda^*\|_{\ell_2} = \frac{1}{\tau^2} \sqrt{d(\lambda^*)} \varepsilon.$$

If $h_j(X)$, $j = 1, \ldots, N$, are i.i.d. $N(0, \tau^2)$, this yields

$$\|\hat{\lambda} - \lambda^*\|_{\ell_2} = \beta_2^2(d(\lambda^*)) \sqrt{d(\lambda^*)} \varepsilon,$$

which is in agreement with the last bound of Theorem 2 in this case. Thus, the presence of $\beta_2(d)$ in the bound has something to do with the nature of the problem, although there might be different and, possibly, much better ways to take into account the geometry of the dictionary.

*Example.* Suppose that

$$\{1, \ldots, N\} := \bigcup_{k=1}^{m} I_k,$$

where $I_k$ are disjoint sets. Suppose that $\phi_1(X), \ldots, \phi_m(X)$ are i.i.d. $N(0, 1)$. Let

$$g_\mu = \sum_{k=1}^{m} \mu_k \phi_k, \qquad \mu = (\mu_1, \ldots, \mu_m) \in \mathbb{R}^m,$$

and define $h_j := \phi_k, j \in I_k$. It is easy to check that, for such a dictionary, the Dantzig selector $\hat{\lambda}$ is a solution of the following problem. First, solve the problem

$$\sum_{j=1}^{m} |\mu_j| \to \min$$

subject to the constraints

$$\max_{1 \leq k \leq m} \left| n^{-1} \sum_{j=1}^{n} (g_\mu(X_j) - Y_j) \phi_k(X_j) \right| \leq \varepsilon.$$

Let $\hat{\mu}$ be its solution. Then, take arbitrary $\hat{\lambda}$ satisfying the conditions

$$\sum_{j \in I_k} \hat{\lambda}_j = \hat{\mu}_k, \qquad \text{sign}(\hat{\lambda}_j) = \text{sign}(\hat{\mu}_k), \qquad j \in I_k, \ k = 1, \ldots, m.$$

Clearly, we have $f_{\hat{\lambda}} = g_{\hat{\mu}}$ and

$$\|f_{\hat{\lambda}} - f_*\|^2_{L_2(\Pi)} = \|g_{\hat{\mu}} - f_*\|^2_{L_2(\Pi)}.$$

If $f_* = g_{\mu^*}$ for some $\mu^* \in \mathbb{R}^m$, then it follows from Theorem 2 that, with a high probability,

$$\|f_{\hat{\lambda}} - f_*\|^2_{L_2(\Pi)} = \|g_{\hat{\mu}} - f_*\|^2_{L_2(\Pi)} \leq Dd(\mu^*)\varepsilon^2$$

(under appropriate further assumptions, say, that $\xi$ is $N(0, \sigma^2)$ and

$$\varepsilon \geq C\sigma\sqrt{\frac{A\log N}{n}}).$$

Of course, in this case, the dictionary is linearly dependent, coefficients of representation $f_\lambda$ are not identifiable and, in addition, such quantities as $\beta_2(J)$ are infinite. But the Dantzig selector is still recovering $f_*$ with the $L_2(\Pi)$-error being of the correct size.

This example shows that there are situations beyond the scope of Theorems 1 and 2 in which the Dantzig selector is a reasonably good estimator of an unknown regression function $f_*$. It might be possible to develop more subtle geometric characteristics of the dictionary than $\beta_2$ which can be used, for instance, when the dictionary can be partitioned into disjoint sets of highly correlated functions with very little correlation between the sets, and to prove sparsity oracle inequalities in terms of such characteristics. However, it is not our goal in this paper to study such situations in detail.

## Appendix: Several exponential bounds

We need the following three lemmas.

**Lemma 3.** *Let* $\eta^{(k)}, \eta_1^{(k)}, \ldots, \eta_n^{(k)}$ *be i.i.d. random variables with* $\mathbb{E}\eta^{(k)} = 0$ *and* $\|\eta^{(k)}\|_{\psi_1} < +\infty$, $k = 1, \ldots, N$. *There exists a numerical constant* $C > 0$ *such that for all* $A \geq 1$ *with probability at least* $1 - N^{-A}$, *for all* $k = 1, \ldots, N$,

$$\left| n^{-1} \sum_{j=1}^n \eta_j^{(k)} \right| \leq C\|\eta^{(k)}\|_{\psi_1} \left( \sqrt{\frac{A\log N}{n}} \vee \frac{A\log N}{n} \right).$$

This is a consequence of a well-known version of Bernstein's inequality (see, for example, Lemma 2.2.11 in van der Vaart and Wellner (1996)).

**Lemma 4.** *There exists a constant* $C > 0$ *such that for all* $A \geq 1$ *with probability at least* $1 - N^{-A}$,

$$\sup_{\|u\|_{\ell_1} \leq 1} |(\Pi_n - \Pi)(|f_u|)| \leq C \max_{1 \leq k \leq N} \|h_k(X)\|_{\psi_1} \left( \sqrt{\frac{A\log N}{n}} \vee \frac{A\log N}{n} \right).$$

**Proof.** Let $R_n(f)$ denote the Rademacher process

$$R_n(f) := n^{-1} \sum_{j=1}^{n} \varepsilon_j f(X_j),$$

$\varepsilon, \varepsilon_j, j = 1, \ldots, n$, being i.i.d. Rademacher random variables independent of $X_1, \ldots, X_n$. For $t > 0$, we use the symmetrization inequality and then the contraction inequality (see Ledoux and Talagrand (1991), page 112) to get

$$\mathbb{E}\exp\Big\{t \sup_{\|u\|_{\ell_1} \leq 1} |(\Pi_n - \Pi)(|f_u|)|\Big\} \leq \mathbb{E}\exp\Big\{2t \sup_{\|u\|_{\ell_1} \leq 1} |R_n(|f_u|)|\Big\}$$

$$\leq \mathbb{E}\exp\Big\{4t \sup_{\|u\|_{\ell_1} \leq 1} |R_n(f_u)|\Big\}.$$

Since the mapping $u \mapsto R_n(f_u)$ is linear, the supremum of $R_n(f_u)$ over the set $\{\|u\|_{\ell_1} \leq 1\}$ is attained at one of its vertices and we get

$$\mathbb{E}\exp\Big\{t \sup_{\|u\|_{\ell_1} \leq 1} |(\Pi_n - \Pi)(|f_u|)|\Big\} \leq \mathbb{E}\exp\Big\{4t \max_{1 \leq k \leq N} |R_n(h_k)|\Big\}$$

$$= N \max_{1 \leq k \leq N} \mathbb{E}[\exp\{4tR_n(h_k)\} \vee \exp\{-4tR_n(h_k)\}]$$

$$\leq 2N \max_{1 \leq k \leq N} \mathbb{E}\exp\{4tR_n(h_k)\}$$

$$\leq 2N \max_{1 \leq k \leq N} \Big(\mathbb{E}\exp\Big\{4\frac{t}{n}\varepsilon h_k(X)\Big\}\Big)^n.$$

To bound the last expectation and to complete the proof, we need only to follow the standard proof of the Bernstein inequality. □

The proof of the next lemma is quite similar.

**Lemma 5.** *There exists a constant $C > 0$ such that for all $A \geq 1$, with probability at least $1 - N^{-A}$,*

$$\sup_{\|u\|_{\ell_1} \leq 1} |(\Pi_n - \Pi)(|f_u|^2)| \leq C \max_{1 \leq k,j \leq N} \|h_k(X)h_j(X)\|_{\psi_1} \left(\sqrt{\frac{A \log N}{n}} \vee \frac{A \log N}{n}\right).$$

Let $J \subset \{1, \ldots, N\}$, $d(J) \leq d \leq \frac{N}{e} - 1$. Define

$$K_J := C_J \cap \{u \in \mathbb{R}^N : \|u\|_{\ell_2} \leq 1\}.$$

**Lemma 6.** *There exists a constant $C > 0$ such that for all $A \geq 1$, with probability at least*

$$1 - 5^{-dA}\binom{N}{\leq d}^{-A},$$

$$\sup_{u \in K_J}|(\Pi_n - \Pi)(|f_u|)| \leq C \sup_{\|u\|_{\ell_2}\leq 1, d(u)\leq d}\|f_u\|_{\psi_1}\left(\sqrt{\frac{Ad\log(N/d)}{n}} \vee \frac{Ad\log(N/d)}{n}\right).$$

**Proof.** The idea of the proof is well known (see, for example, Ledoux and Talagrand (1991), page 421, or, in a context closer to the current paper, Mendelson, Pajor and Tomczak-Jaegermann (2007), Lemma 3.3). Recall Lemma 1 and its notation. Let $u \in K_J$. Lemma 1 implies that

$$K_J \subset 3\operatorname{conv}\left(\bigcup B_I : I \subset \{1,\ldots,N\}, d(I) \leq d\right),$$

where

$$B_I := \left\{(u_i : i \in I) : \sum_{i \in I}|u_i|^2 \leq 1\right\}$$

since $u^{(0)} \in B_{J_0}$, $u^{(1)} \in B_{J_1}$ and

$$\sum_{k \geq 2} u^{(k)} \in \operatorname{conv}\left(\bigcup B_I : I \subset \{1,\ldots,N\}, d(I) \leq d\right).$$

It is easy to see that if $B$ is the unit Euclidean ball in $\mathbb{R}^d$ and $M$ is a 1/2-net of this ball, then

$$B \subset 2\operatorname{conv}(M).$$

A somewhat informal version of the proof of this claim that can easily be made precise is as follows: with $+$ denoting Minkowski sum,

$$B \subset M + \tfrac{1}{2}B \subset \operatorname{conv}(M) + \tfrac{1}{2}B \subset \operatorname{conv}(M) + \tfrac{1}{2}\operatorname{conv}(M) + \tfrac{1}{4}B \subset \cdots$$
$$\subset \operatorname{conv}(M) + \tfrac{1}{2}\operatorname{conv}(M) + \tfrac{1}{4}\operatorname{conv}(M) + \cdots \subset 2\operatorname{conv}(M).$$

Now selecting for each $I$ with $d(I) \leq d$ its minimal 1/2-net $M_I$ yields

$$K_J \subset 6\operatorname{conv}\left(\bigcup M_I : I \subset \{1,\ldots,N\}, d(I) \leq d\right) =: 6\ \operatorname{conv}(\mathcal{M}_d).$$

Therefore, we can repeat the proof of Lemma 4 and reduce the bounding of

$$\sup_{u \in K_J}|(\Pi_n - \Pi)(|f_u|)|$$

to the bounding of

$$\sup_{u \in \mathcal{M}_d} |R_n(f_u)|,$$

with $\text{card}(\mathcal{M}_d)$ playing the role of $N$. It remains to observe that

$$\text{card}(\mathcal{M}_d) \leq 5^d \binom{N}{\leq d},$$

which implies that with some $c > 0$,

$$\log(\text{card}(\mathcal{M}_d)) \leq cd \log \frac{N}{d},$$

and the result follows. □

## Acknowledgements

This work was partially supported by NSF Grant DMS-MSPA-0624841 and NSF Grant CCF-0808863.